\documentclass[11pt,a4paper]{article}
\usepackage[latin1]{inputenc}
\usepackage[hmargin=1.5cm,vmargin=2.5cm]{geometry}
\usepackage{amsthm}
\usepackage{graphicx}
\usepackage{lmodern}
\usepackage{newlfont}
\usepackage[T1]{fontenc}
\usepackage{amsmath}
\usepackage{amsfonts}
\usepackage{amssymb}
\usepackage{graphicx}
\usepackage[english]{babel}
\newtheorem{thm}{Theorem}[section]

\newtheorem{prop}{Proposition}[section]

\newtheorem{dfn}{Definition}[section]

\newtheorem{lm}{Lemma}[section]

\newtheorem{rque}{Remark}[section]
\newtheorem{exple}{Example}[section]
\usepackage{fancyhdr}
\begin{document}

 \date{}
\title{\textbf{Half-Lightlike Submanifold of a Lorentzian Manifolds with a conformal co-screen distribution}}
\author{Issa Allassane Kaboye\footnote{Facult\'e des Sciences et Techniques,Universit\'e de Zinder,BP 656 Zinder( Niger) E-mail: allassanekaboye@yahoo.fr}, Mahamane Mahi Harouna\footnote{Facult\'e des Sciences et Techniques, Universit\'e Dan Dicko Dankoulodo de Maradi( Niger), BP 465 Maradi( Niger),E-mail: hmahi2007@yahoo.fr}\footnote{correspondig author}  and Bazanfar\'e Mahaman\footnote{Facult\'e des Sciences et Techniques, Universit\'e Abdou Moumouni  de Niamey( Niger), BP 10662 Niamey( Niger)
E-mail: bmahaman2007@yahoo.fr}}
\maketitle
\textbf{Abstract.} In this paper, we give the Cartan$^{\prime}$s formula for 
half-lightlike submanifolds of Lorentzian manifolds and use it to show that a screen homothetic  half-lightlike submanifolds of a Lorentzian space form, with a conformal co-screen distribution are locally a lightlike triple product manifolds. 
Then we give a classification theorem for half-lightlike submanifolds of Lorentzian space form with constant screen principal curvatures. These results  extend some results obtained in the case of lightlike hypersurfaces of Lorentzian manifolds (\cite{C.M.J}).

  $MSC2010: 53B30, 53C42, 53C50$

\section{Introduction}
It is well known that the intersection of the normal bundle and the tangent bundle of a lightlike submanifold
of a semi-Riemannian manifold is not trivial(\cite{K.A2}). Thus, one
cannot use, in the usual way, the classical submanifold theory to define any induced
object on a lightlike submanifold. To overcome these difficulties in degenerate geometry, Duggal and Bejancu (\cite{K.A2}) 
introduced a non-degenerate screen distribution (or equivalently a null transversal 
vector bundle) so as to get three factors splitting the ambient tangent bundle. Then, they derived the main induced geometric
objects such as second fundamental forms, shape operators, induced connections, curvature, etc.

The class of lightlike submanifolds of codimension $2$ is composed of two classes by virtue of the 
dimension of its radical distribution (dimension that is either 1 or 2) , named by half-lightlike and coisotropic 
submanifolds (\cite{K.A1}, \cite{K.D}). A codimension 2 lightlike submanifold is called a half-lightlike submanifold if 
dim(Rad(TM)) = 1. For more results about half-lightlike submanifolds, we refer the reader 
to (\cite{D.H.J1}, \cite{D.H.J2}, \cite{K.S}).

We know that the shape operator plays an important role in the study geometry of submanifolds. In the case of half-lightlike submanifolds, there are three shape operators 
($A_N~\stackrel{\ast}{A}_\xi~\mbox{and}~A_L$) and there are interrelations between these geometric objects 
and those of its screen distributions. The shape operators $A_N~\mbox{and}~A_\xi$ are $\Gamma(S(TM))$-valued, but not $A_L$. 
Moreover, the shape operators $A_N~\mbox{and}~A_L$ of a half-lightlike submanifold are not necessarily auto-adjoint, but the 
operator $\stackrel{\ast}{A}_\xi$ of the screen distribution is diagonalizable. 
In \cite{C.M.J}, C.Atindogb\'e, M.M. Harouna and J.Tossa prove the so-called Cartan$^{,}$s fundamental formula for lightlike 
hypersurfaces and use it to show that a screen conformal lightlike hypersurface of a Lorentzian Euclidean space is locally a 
lightlike triple product manifold. They also give a classification theorem for lightlike hypersurfaces of Lorentzian Euclidean 
space with constant screen principal curvatures. In this paper, we generalize those studies on half-lightlike submanifold 
case. The article is organised as follows: in section $2$ we brief  basic informations on half-lightlike submanifolds following
closely the approach in (\cite{K.A2}). In section 3,  we  prove the  Cartan's  formula for half-lightlike submanifold. 
 In section 4, we use this formula to show that a screen homothetic  half-lightlike submanifolds 
of a Lorentzian space form, with a conformal co-screen distribution are locally a lightlike triple product manifolds. 
Then we give a classification theorem for half-lightlike submanifolds of Lorentzian space form with constant screen 
principal curvatures.
\section{Preliminaries on half-lightlike submanifolds}
In this section, we present basic notions on  differential geometry of half-lightlike submanifold manifolds. A full discussion of 
this content  can be found in (\cite[]{K.A1, K.A2, K.D, K.S}).

Let $(\overline{M},\overline{g})$ be a $(m+2)$-dimensional semi-Riemannian manifold of index $q\geq 1$
and $(M,g)$ a lightlike submanifold of codimension 2 of $\overline{M}$. We say that $M$ is a half-lightlike submanifold if the dimension of the radical distribution $Rad(TM)$ is one. It is well known that the radical distribution is given by $Rad(TM) = TM \cap TM^\bot$, where $TM^\bot$ is called the normal bundle of $M$ in $\overline{M}$. Thus there exist two non-degenerate complementary distributions $S(TM)$ and
$S(TM^\bot)$ of $Rad(TM)$ in $TM$ and $TM^\bot$ respectively \cite{K.A2}, which are called the screen and co-screen distribution respectively on $M$ . Thus we have
\begin{equation}\label{Prelimeq1}
TM=\mbox{Rad}(TM)\underset{\bot}{\oplus} S(TM)
\end{equation}
and
\begin{equation}\label{Prelimeq2}
TM^{\bot}=\mbox{Rad}(TM)\underset{\bot}{\oplus} S(TM^{\bot})
\end{equation}
where $\underset{\bot}{\oplus}$ denotes the orthogonal direct sum.
\par Since  $S(TM)$ is non-degenerate, we have the decomposition
\begin{equation}\label{Prelimeq3}
T\overline{M}_{\mid_{M}}=S(TM)\underset{\bot}{\oplus} S(TM)^{\bot}
\end{equation}
Thus, it is easy to see that $S(TM^\bot)$ is a subbundle of $S(TM)^\bot$. As $S(TM^\bot)$ is a non-degenerate
subbundle of $S(TM)^\bot$, we have the decomposition
\begin{equation}\label{Prelimeq4}
S(TM)^{\bot}=S(TM^\bot)^\bot\underset{\bot}{\oplus} S(TM^{\bot})
\end{equation}
It is obvious to see that $S(TM^\bot)^\bot$ is also a non-degenerate distribution and $Rad(TM)$ is
a subbundle of $S(TM^\bot)^\bot$. Choose $L\in \Gamma(S(TM^\bot))$ such that $\overline{g}(L,L)=1$. For any null section $\xi \in (Rad(TM))$, there exists a unique null vector field $N\in (S(TM^\bot)^\bot)$ satisfying
\begin{eqnarray}\label{Prelimeq5}
\overline{g}(\xi,N) =1, \overline{g}(N,N)= \overline{g}(N,X)= \overline{g}(N,L)= 0, \forall X\in (S(TM)).
\end{eqnarray}
Denote by $ltr(TM)$ the vector subbundle of $S(TM^\bot)^\bot$ locally spanned by $N$. Then we have:
\begin{eqnarray}
 S(TM^\bot)^\bot = Rad(TM)\oplus ltr(TM).
\end{eqnarray}
 Let $tr(TM) = S(TM^\bot)\underset{\bot}{\oplus} ltr(TM)$. We call $N$, $ltr(TM)$ and $tr(TM)$ the lightlike transversal vector field, lightlike transversal vector bundle and transversal vector bundle of M with respect to the chosen screen distribution $S(TM)$ respectively. Then $T\overline{M}$ is decomposed as following
\begin{eqnarray}\label{Prelimeq6}
T\overline{M}_{\mid_{M}} &=& TM\oplus tr(TM)= Rad(TM)\oplus tr(TM)\underset{\bot}{\oplus}S(TM)
\cr&=& Rad(TM)\oplus ltr(TM)\underset{\bot}{\oplus}S(TM)\underset{\bot}{\oplus} S(TM^\bot).
\end{eqnarray}
Let $P$ be the projection morphism of $TM$ on $S(TM)$ with respect to the decomposition (\ref{Prelimeq1}). For any $X,Y\in\Gamma(TM)$, $N\in \Gamma(ltr(TM))$, $\xi\in\Gamma(Rad(TM))$ and $L\in\Gamma(S(TM)^\bot)$, the Gauss and Weingarten formulas of $M$ and $S(TM)$ are given by
\begin{eqnarray}
\overline{\nabla}_X Y &=& \nabla_X Y+B(X,Y)N+D(X,Y)L,\label{GW1}\\
\overline{\nabla}_X N &=& -A_N X+\tau(X)N+\rho(X)L,\label{GW2}\\
\overline{\nabla}_X L &=& -A_L X+\phi(X)N,\label{GW3}\\
\nabla_X PY &=&\stackrel{\ast} {\nabla}_X PY+C(X,PY)\xi,\label{GW4}\\
\nabla_X \xi &=& -\stackrel{\ast} {A}_\xi X-\tau(X)\xi,\label{GW5}
\end{eqnarray}
respectively, where $\nabla ~\mbox{and}~\stackrel{\ast}{\nabla}$ are induced connection on $TM$ and $S(TM)$ respectively, $B$ and $D$ are called local second fundamental forms of $M$, $C$ is called the local second fundamental form of $S(TM)$; $A_N,~\stackrel{\ast}{A}_\xi~\mbox{and}~ A_L$ are linear shape operators on $TM$ of lightlike transversal bundle, radical bundle and screen transversal bundle respectively. $\tau,~ \rho~\mbox{and}~ \phi$ are 1-forms on $M$.
\begin{eqnarray}\label{secondfond}
h(X,Y) = B(X,Y)N + D(X,Y)L
\end{eqnarray}
is the second fundamental form tensor of $M$. Since the connection $\overline{\nabla}$ on $\overline{M}$ is torsion-free, the induced connection $\nabla$ on $M$ is also torsion-free, and then $B$ and $D$ are symmetric tensors on $\Gamma(TM)$. But $\nabla$ is not metric, since using (\ref{GW1}), we have
\begin{eqnarray}\label{metric}
(\nabla_X g)(Y ,Z)=B(X,Y )\eta(Z)+B(X,Z)\eta(Y ),
\end{eqnarray}
for all $X,Y \in\Gamma(TM)$, where $\eta$ is a 1-form on $TM$ such that $\eta(X) = g(X,N)$,~~$\forall X\in\Gamma(TM)$. It is easy to verify that the induced connection $\stackrel{\ast} {\nabla}$ on $S(TM)$ is metric. Since $B(X,Y)=\overline{g}(\overline{\nabla}_X Y ,\xi)~\mbox{and}~D(X,Y)=\overline{g}(\overline{\nabla}_X Y ,L)$, it is obvious that $B~\mbox{and}~D$ are independent of the choice of a screen distribution. However, we note that both $B$ and $\tau$ depend on the section $\xi\in\Gamma(Rad(TM))$. In the case we take $\xi^*=\alpha\xi$, it follows that  $N^{*}=\frac{1}{\alpha}N$. Hence we obtain 
\begin{eqnarray}\label{Prelimeq7}
B^*=\alpha B~~\mbox{and}~~\tau(X)=\tau^*(X)+X(\log\alpha)~~\mbox{for any}~~X\in\Gamma(TM)
\end{eqnarray}
The above three local second fundamental forms of $M$ and $S(TM)$ are related to their shape operators by
\begin{eqnarray}
g(\stackrel{\ast}{A}_\xi X,Y)&=& B(X,Y),~~~~~~~~~~~~~~~~~~~~~~\overline{g}(\stackrel{\ast}{A}_\xi X,N)=0\label{shape01}\\
g(A_L X,Y)&=& D(X,Y)+\phi(X)\eta(Y),~~~~~\overline{g}(A_L X,N)=\rho(X),\label{shape02}\\
g(A_N X,PY)&=& C(X,PY),~~~~~~~~~~~~~~~~~~~\overline{g}(A_N X,N)=0,\label{shape03}
\end{eqnarray} 
for any $X,Y\in\Gamma(TM)$.\\
From (\ref{shape01}),~(\ref{shape02}),~(\ref{shape03}), we see that $B$ and $D$ satisfy
\begin{eqnarray}\label{secondformfondprop}
B(X, \xi) = 0,~~D(X, \xi) = -\phi(X),~\forall X \in\Gamma(TM),
\end{eqnarray}
$\stackrel{\ast}{A}_\xi~\mbox{and}~A_N$ are $S(TM)$-valued, and $\stackrel{\ast}{A}_\xi$ is self-adjoint on $TM$ such that
\begin{eqnarray}
\stackrel{\ast}{A}_\xi\xi=0
\end{eqnarray}
Denote by $\overline{R},~R~\mbox{and}~\stackrel{\ast}{R}$ the curvature tensors of connection
$\overline{\nabla}~\mbox{on}~\overline{M}$, the induced connection $\nabla$ on $M$ and the induced connection $\stackrel{\ast}{\nabla}$ on $S(TM)$, respectively. Using the Gauss-Weingarten formulas for $M$ and $S(TM)$, we obtain the Gauss-Codazzi equations for $M$ and $S(TM)$:
\begin{eqnarray}\label{GausCoda1}
 \overline{R}(X,Y)Z&=& R(X,Y)Z+B(X,Z)A_N Y-B(Y,Z)A_N X
 \cr&  & +D(X,Z)A_L Y-D(Y,Z)A_L X
 \cr&  & +\{(\nabla_X B)(Y,Z)-(\nabla_Y B)(X,Z)
 \cr&  & + B(Y,Z)\tau(X)-B(X,Z)\tau(Y)
 \cr&  & +D(Y,Z)\phi(X)-D(X,Z)\phi(Y)\}N
 \cr&  & +\{(\nabla_X D)(Y,Z)-(\nabla_Y D)(X,Z)+B(Y,Z)\rho(X)
 \cr&  & -B(X,Z)\rho(Y)-D(Y,Z)+D(X,Z)\}L,
\end{eqnarray} 
\begin{eqnarray}\label{GausCoda2}
\overline{R}(X,Y)N&=& -\nabla_X( A_N Y)+\nabla_Y( A_N X)+A_N[X,Y]
\cr&  & +\tau(X)A_N Y-\tau(Y)A_N X+\rho(X)A_L Y-\rho(Y)A_L X
\cr&  & +\{B(A_N X,Y)-B(A_N Y,X)+2d\tau(X,Y)+
\cr&  & +\phi(X)\rho(Y)-\phi(Y)\rho(X)\}N
\cr&  & +\{D(A_N X,Y)-D(A_N Y,X)+2d\rho(X,Y)+\rho(X)\tau(Y)-\rho(Y)\tau(X)\}L,
\end{eqnarray}
\begin{eqnarray}\label{GausCoda3}
\overline{R}(X,Y)L&=& -\nabla_X( A_L Y)+\nabla_Y( A_L X)+A_L[X,Y]
\cr&  &+ \phi(X)A_N Y-\phi(Y)A_N X
\cr&  &+ \{B(A_L X,Y)-B(A_L Y,X)+2d\phi(X,Y)
\cr&  &+ \tau(X)\phi(Y)-\tau(Y)\phi(X)\}N
\cr&  &+\{D(A_L X,Y)-D(A_L Y,X)+\rho(X)\phi(Y)-\rho(Y)\phi(X)\}L,
\end{eqnarray}
\begin{eqnarray}\label{GausCoda4}
R(X,Y)PZ &=& \stackrel{\ast}{R}(X,Y)PZ+C(X,PZ)\stackrel{\ast}{A}_\xi Y-C(Y,PZ)\stackrel{\ast}{A}_\xi X
\cr&  &+\{(\nabla_X C)(Y,PZ)-(\nabla_Y C)(X,PZ)
\cr&  &+C(X,PZ)\tau(Y)-C(Y,PZ)\tau(X)\}\xi,
\end{eqnarray} 
\begin{eqnarray}\label{GausCoda5}
R(X,Y)\xi &=& -\stackrel{\ast}{\nabla}_X(\stackrel{\ast}{A}_\xi Y)+\stackrel{\ast}{\nabla}_Y(\stackrel{\ast}{A}_\xi X)+\stackrel{\ast}{A}_\xi [X,Y]+\tau(Y)\stackrel{\ast}{A}_\xi X
\cr&  & -\tau(X)\stackrel{\ast}{A}_\xi Y+\{C(Y,\stackrel{\ast}{A}_\xi X)-C(X,\stackrel{\ast}{A}_\xi Y)-2d\tau(X,Y)\}\xi.
\end{eqnarray}
A semi-Riemannian manifold $\overline{M}$ of constant curvature $k$ is called a semi-Riemannian
space form and denote it by $\overline{M}(k)$. The curvature tensor $\overline{R}$ of $\overline{M}(k)$ is given by
\begin{eqnarray}\label{courbeq1}
\overline{R}(\overline{X},\overline{Y})\overline{Z} = k\{\overline{g}(\overline{Y},\overline{Z})\overline{X} -\overline{g}(\overline{X},\overline{Z})\overline{Y} \},~\forall \overline{X},\overline{Y},\overline{Z}\in\Gamma(T\overline{M}).
\end{eqnarray}
It is then obvious that $\overline{R}(X,Y)Z$ is tangent to $M$ for all $X,Y,Z\in\Gamma(TM)$ , so from (\ref{GausCoda1}), we get
\begin{eqnarray}
\overline{R}(X,Y)Z = R(X,Y)Z + B(X,Z)A_N Y - B(Y,Z)A_N X
+ D(X,Z)A_L Y - D(Y,Z)A_L X,~\forall X,Y,Z \in\Gamma(TM).
\end{eqnarray}


 The Ricci curvature tensor of $\overline{M}$, denoted by $\overline{R}$ic is defined by
\begin{eqnarray}\label{ricdengeneq1}
\overline{R}ic=trace\{Z\rightarrow \overline{R}(Z,X)Y\},~~\mbox{for any}~~ X,Y\in\Gamma(T\overline{M}).
\end{eqnarray}
Locally, $\overline{R}$ic is given by
\begin{eqnarray}\label{ricdengeneq2}
\overline{R}ic=\sum_{i=1}^{m+2}\varepsilon_{i}\overline{g}(\overline{R}(E_i,X)Y,E_i),
\end{eqnarray}
where $\{E_1,\cdots, E_{m+2}\}$ is an orthonormal frame field of $T\overline{M}$ and $\varepsilon_{i}=\pm 1$ denote the causal caracter of respective vector field $E_i$.
\par Consider the induced quasi-orthonormal frame field $\{\xi, W_a\}$ on $M$, where $Rad(TM)=span\{\xi\}$ and $S(TM)=span\{W_a\}$ and let $E=\{\xi; W_a, L, N\}$ be the corresponding frame field on $\overline{M}$. Then, by using (\ref{ricdengeneq1}), we obtain
\begin{eqnarray}\label{ricdengeneq3}
\overline{R}ic(X,Y)=\sum_{a=1}^{m-1}\varepsilon_{a}\overline{g}(\overline{R}(W_a,X)Y,W_a)+\overline{g}(\overline{R}(\xi ,X)Y,N)+\overline{g}(\overline{R}(L ,X)Y,L)+\overline{g}(\overline{R}(N ,X)Y,\xi).
\end{eqnarray}
Let $Ric$ denote the induced Ricci tensor on $M$ given by:
\begin{eqnarray}\label{ricdengeneq4}
Ric(X,Y)=trace\{Z\rightarrow R(Z,X)Y\}~~\mbox{for any}~~X,Y\in\Gamma(TM). 
\end{eqnarray}
Using the quasi-orthonormal frame $\{\xi, W_a\}$ on $M$, we obtain:
\begin{eqnarray}\label{ricdengeneq5}
Ric(X,Y)=\sum_{a=1}^{m-1}\varepsilon_{a}g(R(W_a,X)Y,W_a)+\overline{g}(R(\xi ,X)Y,N).
\end{eqnarray}
Substituting the Gauss-Codazzi equations (\ref{GausCoda1}) in (\ref{ricdengeneq3}) and using relation (\ref{shape01}) and (\ref{shape03}), we obtain
\begin{eqnarray}\label{ricdengeneq6}
Ric(X)=\overline{R}ic(X,Y)+B(X,Y)tr(A_N)+D(X,Y)tr(A_L)-g(A_{N}X,\stackrel{\ast}{A}_{\xi}Y)-g(A_{L}X,A_{L}Y)\nonumber\\
+\rho(X)\varphi(Y)-\overline{g}(\overline{R}(\xi ,Y)X,N)-\overline{g}(\overline{R}(L ,Y)X,L)~~\mbox{for any}~~X,Y\in\Gamma(TM).
\end{eqnarray}
This shows that $Ric$ is not symmetric. Using (\ref{ricdengeneq6}) and the first Bianchi's identity, we obtain:
\begin{eqnarray}\label{ricdengeneq7}
Ric(X,Y)-Ric(Y,X)=g(\stackrel{\ast}{A}_{\xi}X,A_{N}Y)-g(\stackrel{\ast}{A}_{\xi}Y,A_{N}X)+\rho(X)\varphi(Y)-\rho(Y)\varphi(X)-\overline{g}(\overline{R}(X,Y)\xi,N).
\end{eqnarray}
From this relation and the relations (\ref{GausCoda1}) and (\ref{GausCoda5}), we have
\begin{eqnarray}\label{ricdengeneq8}
Ric(X,Y)-Ric(Y,X)=2d\tau(X,Y).
\end{eqnarray}
\begin{thm}(\cite{K.S})\label{theorem1form}
\textsl{Let $(M,g,S(TM),S(TM^\bot))$ be a half-lightlike submanifold of a semi-Riemannian manifold $(\overline{M},\overline{g})$. Then, the induced Ricci curvature on $M$ is symmetric if and only if each 1-form $\tau$ is closed, i.e $d\tau=0$ on any neighborhood $\mathcal{U}\subset M$. }
\end{thm}
\begin{rque}\label{rque3}
\textsl{If the induced Ricci tensor $Ric$, on $M$ is symmetric, the 1-form $\tau$ is closed by 
theorem\ref{theorem1form}. Thus, by Poincar\'e lemma, there existe a smooth function $f$ on $\mathcal{U}\subset M$ such that $\tau=df$. Consequently, we get $\tau(X)=X.(f)~~\mbox{for any}~~X\in\Gamma(TM_{\mid_{\mathcal{U}}})$. If we take $\xi^{*}=\alpha\xi$, then by setting $\alpha=\exp(f)$, we get $\xi^{*}=0~~\mbox{for any}~~X\in\Gamma(TM_{\mid_{\mathcal{U}}})$. We call the pair~ $\{\xi,N\}$~ on ~$\mathcal{U}$~ such that the corresponding 1-form $\tau$ vanishes the canonical null pair of $M$.}
\end{rque},

\begin{dfn}\textsl{
A half-lightlike submanifold $(M,g,S(TM),S(TM^\bot))$ of a semi-Riemannian manifold is screen locally (resp., globally) 
conformal if on any coordinate neighborhood $\mathcal{U}$ (resp., $\mathcal{U}=M$) there exists a non-zero smooth 
function $\varphi$ such that for any null vector field $\xi\in\Gamma(TM^\bot)$, the relation
\begin{equation}\label{conf1}
A_N X=\varphi\stackrel{\ast}{A}_\xi X,~~\mbox{for any}~~X\in\Gamma(TM_{\vert\mathcal{U}}), 
\end{equation}
holds between the shape operators $A_N~~\mbox{and}~~\stackrel{\ast}{A}_\xi$ of $M~~\mbox{and}~~S(TM)$, respectively.}

 \textsl{In particular, if $\varphi$ is a non-zero constant, $M$ is called screen homothetic.}
 \end{dfn}
 It is easy to see from (\ref{shape01}) and (\ref{shape03}) that a half-lightlike submanifold $M$ is screen conformal if and 
 only if the second fundamental forms $B$ and $C$ satisfy
\begin{equation}
C(X,PY)=\varphi B(X,Y),~~\mbox{for any}~~X,Y\in\Gamma(TM_{\mid_{\mathcal{U}}})
\end{equation}

\section{Cartan$^{\prime}$s formula for half-lightlike submanifolds}
In this section, we first consider a half-lightlike submanifold $(M,g,S(TM),S(TM^\bot))$ of a semi-Riemannian manifold $(\overline{M}(k),\overline{g})$
of constant curvature $k$. As a direct application of the Gauss-Codazzi equations for $M~~\mbox{and}~~S(TM)$, we have this proposition.
\begin{prop}\label{prop1}
\textsl{Let $(\overline{M}(k),\overline{g})$ be a semi-Riemannian manifold of constant curvature $k$ and $(M,g,S(TM),S(TM^\bot))$ be a
half-lightlike submanifold of $\overline{M}(k)$. Denote by $R$ the curvature tensor of the induced connection $\nabla$ on $M$ by the Levi-civita connection $\overline{\nabla}$ on $\overline{M}(k)$. For any $X,Y,Z \in\Gamma(TM)$, we have:
\begin{enumerate}
\item $R(X,Y)Z=k\{g(Y,Z)X-g(X,Z)Y\}+B(Y,Z)A_N X-B(X,Z)A_N Y+D(Y,Z)A_L X-D(X,Z)A_L Y$
\item $(\nabla_X B)(Y,Z)-(\nabla_Y B)(X,Z)=B(X,Z)\tau(Y)-B(Y,Z)\tau(X)+D(X,Z)\phi(Y)-D(Y,Z)\phi(X)$
\item $(\nabla_X D)(Y,Z)-(\nabla_Y D)(X,Z)=B(X,Z)\rho(Y)-B(Y,Z)\rho(X)$
\item $B(A_N Y,X)-B(A_N X,Y)=\phi(X)\rho(Y)-\phi(Y)\rho(X)+2d\tau(X,Y)$
\item $(\nabla_X A_N)Y-(\nabla_Y A_N)X+k[\eta(Y)X-\eta(X)Y]=\tau(X)A_N Y-\tau(Y)A_N X+\rho(X)A_L Y-\rho(Y)A_L X$
\item $D(A_N Y,X)-D(A_N X,Y)=\rho(X)\tau(Y)-\rho(Y)\tau(X)+2d\rho(X,Y)$
\item $(\nabla_X A_L)Y-(\nabla_Y A_L)X=\phi(X)A_N Y-\phi(Y)A_N X$
\item $B(A_L Y,X)-B(A_L X,Y)=\tau(X)\phi(Y)-\tau(Y)\phi(X)+2d\phi(X,Y)$
\item $D(A_L Y,X)-D(A_L X,Y)=\rho(X)\phi(Y)-\rho(Y)\phi(X)$
\item $(\nabla_X \stackrel{\ast}{A}_\xi)Y-(\nabla_Y \stackrel{\ast}{A}_\xi)X=D(X,\xi)A_L Y-D(Y,\xi)A_L X+\tau(Y)\stackrel{\ast}{A}_\xi X-\tau(X)\stackrel{\ast}{A}_\xi Y-2d\tau(X,Y)\xi$
\end{enumerate}
}
\end{prop}
\begin{dfn}
\textsl{A lightlike submanifold $M$ is said to be irrotational if~~$\overline{\nabla}_{X}\xi\in\Gamma(TM)~~\mbox{for any}~~ X\in\Gamma(TM)$, where $\xi\in\Gamma(Rad(TM))$}.
\par \textsl{For a half-lightlike $M$, since $B(X,\xi)=0$, the above definition is equivalent to $D(X,\xi)=\phi(X)=0~~\mbox{for any}~~X\in\Gamma(TM)$}
\end{dfn}
\begin{rque}\label{rque1}
\textsl{From relations (\ref{GW1}) and (\ref{eqqconfcscren}), we get that any half-lightlike submanifold  $(M,g,S(TM),S(TM^\bot))$ of a semi-Riemannian manifold $(\overline{M},\overline{g})$, with a conformal co-screen distribution is an irrotational submanifold, but the opposite way is wrong.}
\end{rque}
In the sequel, we consider a half-lightlike submanifold $M$ of an (m+2)-dimensional Lorentzian manifold $(\overline{M}(k),\overline{g})$ of constant curvature $k$. For this class of screen conformal and irrotational half-lightlike submanifold $M$, the screen distribution $S(TM)$ is Riemannian (Proposition4.2.2 of \cite{K.S}),  integrable ( Theorem4.4.4 of \cite{K.S}) and the induced Ricci tensor on $M$ is symmetric ( Corollaire4.4.11 of  \cite{K.S} ).  Then, according to remark\ref{rque3}, there exists a pair $\{\xi,N\}$ on $\mathcal{U}$ satisfying (\ref{Prelimeq5}) such that the corresponding $1$-forme $\tau$ from (\ref{GW2}) vanishes. Since $\xi$ is an eigenvector field of $\stackrel{\ast}{A}_\xi$ corresponding to the eigenvalue $0$ and $\stackrel{\ast}{A}_\xi$ is $\Gamma(S(TM))$-valued real symmetric, $\stackrel{\ast}{A}_\xi$ has ~$m-1$ ~orthonormal eigenvector fields in $S(TM)$ and is diagonalizable. Consider a frame field of eigenvectors $\{\xi,E_1,\cdots, E_{m-1}\}$ of $\stackrel{\ast}{A}_\xi$ such that $\{E_1,\cdots, E_{m-1}\}$ is an orthonormal frame field of $S(TM)$. Then, $\stackrel{\ast}{A}_\xi E_i=\lambda_i$, $1\leq i\leq m-1$. We call the eigenvalues $\lambda_i$ the screen principal curvatures for all $i$. 
\par In the following, we assume that all screen principal curvatures are constant along $S(TM)~\mbox{and}~\tau=0$. Consider the following distribution on $M$: 
\begin{eqnarray}
T_\lambda = \{X\in\Gamma(S(TM)):~~\stackrel{\ast}{A}_\xi X=\lambda X \}
\end{eqnarray}
 \begin{lm}\label{lm1}
\textsl{Let $(M,g,S(TM),S(TM^\bot))$ be an irrotational half-lightlike submanifold of a semi-Riemannian manifold $(\overline{M}(k),\overline{g})$. If we assume that the 1-forme $\tau$ from (\ref{GW2}) is identically null, then for any $ X\in\Gamma(TM)$, we have:
\begin{enumerate}
\item[(i)] $(\nabla_X \stackrel{\ast}{A}_\xi)Y=(\nabla_Y \stackrel{\ast}{A}_\xi)X~~ \mbox{for any}~~Y\in\Gamma(TM)$
\item[(ii)] $(\stackrel{\ast}{\nabla}_X \stackrel{\ast}{A}_\xi)Y=(\stackrel{\ast}{\nabla}_Y \stackrel{\ast}{A}_\xi)X~~ \mbox{for any}~~Y\in\Gamma(TM)$
\item[(iii)] $\nabla_X \stackrel{\ast}{A}_\xi$ is symmetric with respect to $g$, i.e
\begin{eqnarray*}
g((\nabla_X \stackrel{\ast}{A}_\xi)Y,Z)=g(Y,(\nabla_X \stackrel{\ast}{A}_\xi)Z)~~\forall~~Y,Z\in\Gamma(S(TM))
\end{eqnarray*} 
\item[(iv)]  $\mbox{for any}~~ Y,Z\in\Gamma(S(TM))$
\begin{eqnarray*}
g((\nabla_X \stackrel{\ast}{A}_\xi)Y,Z)&=& g(Y,(\nabla_Z \stackrel{\ast}{A}_\xi)X)= g((\nabla_Z \stackrel{\ast}{A}_\xi)Y,X)
\end{eqnarray*}
\item[(v)]   $\mbox{for any}~~Y\in\Gamma(T_\lambda), Z\in\Gamma(T_\mu)$ 
\begin{eqnarray*}
g((\nabla_X \stackrel{\ast}{A}_\xi)Y,Z)=(\lambda-\mu)g(\nabla_X Y,Z)
\end{eqnarray*}
\end{enumerate}
}
\end{lm}
\textbf{Proof}~~~~\par
Since the  and $\tau=0$,  $(10)$ of proposition\ref{prop1} give $(i)$ .
\begin{eqnarray}\label{20}
(\stackrel{\ast}{\nabla}_X \stackrel{\ast}{A}_\xi)Y &=& \stackrel{\ast}{\nabla}_X (\stackrel{\ast}{A}_\xi Y)-\stackrel{\ast}{A}_\xi\stackrel{\ast}{\nabla}_X Y
\cr&=& \nabla_X\stackrel{\ast}{A}_\xi Y-C(X,\stackrel{\ast}{A}_\xi Y)\xi-\stackrel{\ast}{A}_\xi\nabla_X Y +C(X,Y)\stackrel{\ast}{A}_\xi\xi
\cr&=& \nabla_X\stackrel{\ast}{A}_\xi Y-C(X,\stackrel{\ast}{A}_\xi Y)\xi-\stackrel{\ast}{A}_\xi\nabla_X Y 
\cr&=& (\nabla_X \stackrel{\ast}{A}_\xi)Y-C(X,\stackrel{\ast}{A}_\xi Y)\xi
\cr&=& (\nabla_Y \stackrel{\ast}{A}_\xi)X-C(X,\stackrel{\ast}{A}_\xi Y)\xi
\cr&=& (\stackrel{\ast}{\nabla}_Y\stackrel{\ast}{A}_\xi)X+C(Y,\stackrel{\ast}{A}_\xi X)\xi-C(X,\stackrel{\ast}{A}_\xi Y)\xi
\end{eqnarray}
From (\ref{courbeq1}), it is obvious that $R(X,Y)\xi=0$; and then equation (\ref{GausCoda5}) give 
\begin{eqnarray}\label{21}
g(R(X,Y)\xi, N) = C(Y,\stackrel{\ast}{A}_\xi X) -C(X,\stackrel{\ast}{A}_\xi Y)=0
\end{eqnarray}
 From equations (\ref{20}) and (\ref{21})  we get $(ii)$.
\begin{eqnarray}
g((\nabla_X \stackrel{\ast}{A}_\xi)Y,Z)&=& g(\nabla_X (\stackrel{\ast}{A}_\xi Y),Z)-g(\stackrel{\ast}{A}_\xi \nabla_X Y,Z)
\cr&=& g(\nabla_X (\stackrel{\ast}{A}_\xi Y),Z)-g(\nabla_X Y,\stackrel{\ast}{A}_\xi Z)
\cr&=& -(\nabla_X g)(\stackrel{\ast}{A}_\xi Y,Z)+X.g(\stackrel{\ast}{A}_\xi Y,Z)-g(\stackrel{\ast}{A}_\xi Y,\nabla_X Z)
\cr&  & +(\nabla_X g)(Y,\stackrel{\ast}{A}_\xi Z)-X.g(Y,\stackrel{\ast}{A}_\xi Z)+g(Y,(\nabla_X (\stackrel{\ast}{A}_\xi Z))
\cr&=& -g( Y, \stackrel{\ast}{A}_\xi\nabla_X Z)+g(Y,(\nabla_X (\stackrel{\ast}{A}_\xi Z))
\cr&=& g(Y,(\nabla_X \stackrel{\ast}{A}_\xi)Z)
\end{eqnarray}
That give as $(iii)$.
\begin{eqnarray}\label{22}
g((\nabla_X \stackrel{\ast}{A}_\xi)Y,Z)&=& g(Y,(\nabla_X \stackrel{\ast}{A}_\xi)Z)
\cr&=& g(Y,(\nabla_Z \stackrel{\ast}{A}_\xi)X)
\cr&=& g((\nabla_Z \stackrel{\ast}{A}_\xi)Y,X)
\end{eqnarray}
From $(i)$, $(iii)$ and (\ref{22}), we get $(iv)$
\begin{eqnarray}
g((\nabla_X \stackrel{\ast}{A}_\xi)Y,Z)&=& g(\nabla_X (\stackrel{\ast}{A}_\xi Y),Z)-g(\stackrel{\ast}{A}_\xi\nabla_X Y,Z)
\cr&=& \lambda g(\nabla_X Y,Z)-g(\nabla_X Y,\stackrel{\ast}{A}_\xi Z)
\cr&=& \lambda g(\nabla_X Y,Z)-\mu g(\nabla_X Y,Z)
\cr&=& (\lambda-\mu)g(\nabla_X Y,Z)
\end{eqnarray}
That give $(v)$.
\begin{lm}\label{lm2}
\textsl{Let $(M,g,S(TM),S(TM^\bot))$ be an irrotational half-lightlike submanifold of a semi-Riemannian manifold $(\overline{M}(k),\overline{g})$. If we assume that the 1-forme $\tau$ from (\ref{GW2}) is identically null, then
\begin{enumerate}
\item $\stackrel{\ast}{\nabla}_X Y\in\Gamma(T_\lambda)~~\mbox{if}~~X,Y\in\Gamma(T_\lambda)$ 
\item $\nabla_X Y \bot T_\lambda;~~~~\nabla_Y X \bot T_\mu~~\mbox{if}~~ X\in\Gamma(T_\lambda),~~Y\in\Gamma(T_\mu),~~\lambda\neq\mu.$
\end{enumerate}
}
\end{lm}
\textbf{Proof}\\
Let $Z\in\Gamma(TM)~\mbox{and}~X,Y\in\Gamma(T_\lambda)$. From $(ii)$ and $(iv)$ of lemma\ref{lm1}, we have:
\begin{eqnarray}
g(\stackrel{\ast}{A}_\xi\stackrel{\ast}{\nabla}_X Y,Z)&=& g(\stackrel{\ast}{\nabla}_X \stackrel{\ast}{A}_\xi Y,Z)-g((\stackrel{\ast}{\nabla}_X \stackrel{\ast}{A}_\xi)Y,Z)
\cr&=& \lambda g(\stackrel{\ast}{\nabla}_X Y,Z)-g(Y,(\stackrel{\ast}{\nabla}_Z \stackrel{\ast}{A}_\xi)X)
\cr&=& \lambda g(\stackrel{\ast}{\nabla}_X Y,Z)-g(Y,\stackrel{\ast}{\nabla}_Z \stackrel{\ast}{A}_\xi X)+g(\stackrel{\ast}{\nabla}_Z X,\stackrel{\ast}{A}_\xi Y)
\cr&=&  \lambda g(\stackrel{\ast}{\nabla}_X Y,Z)-\lambda g(Y,\stackrel{\ast}{\nabla}_Z X)+\lambda g(\stackrel{\ast}{\nabla}_Z X,Y)
\cr&=& \lambda g(\stackrel{\ast}{\nabla}_X Y,Z),
\end{eqnarray}
then 
\begin{eqnarray}
g(\stackrel{\ast}{A}_\xi\stackrel{\ast}{\nabla}_X Y-\lambda\stackrel{\ast}{\nabla}_X Y,Z)=0
\end{eqnarray}
Since $Z$ is an arbitrary vector field tangent to $M$, we have:
\begin{eqnarray}
\stackrel{\ast}{A}_\xi\stackrel{\ast}{\nabla}_X Y-\lambda\stackrel{\ast}{\nabla}_X Y=\alpha\xi,~\mbox{where $\alpha$ is a smooth function on $M$.}
\end{eqnarray}
By the fact that
\begin{eqnarray}
\overline{g}(\stackrel{\ast}{A}_\xi\stackrel{\ast}{\nabla}_X Y-\lambda\stackrel{\ast}{\nabla}_X Y,N)=0,
\end{eqnarray}
we conclude that $\alpha=0$, and then we obtain $\stackrel{\ast}{A}_\xi\stackrel{\ast}{\nabla}_X Y=\lambda\stackrel{\ast}{\nabla}_X Y$; and then, we conclude that  $\stackrel{\ast}{\nabla}_X Y\in\Gamma(T_\lambda)$.\\
Now, let $X,Z\in\Gamma(T_\lambda)~\mbox{and}~Y\in\Gamma(T_\mu)$. From $(v)$ of  lemma\ref{lm1}, we have:
\begin{eqnarray}\label{lm21}
g((\nabla_X\stackrel{\ast}{A}_\xi) Y,Z)=(\mu - \lambda)g(\nabla_X Y,Z)
\end{eqnarray}
Otherwise, from $(iv)$ of lemma\ref{lm1} we have
\begin{eqnarray}\label{lm22}
g((\nabla_X\stackrel{\ast}{A}_\xi)Y,Z)=g((\nabla_Z\stackrel{\ast}{A}_\xi)X,Y)=-(\mu-\lambda)g(\nabla_Z X, Y)=-(\mu - \lambda)g(\stackrel{\ast}{\nabla}_Z X,Y).
\end{eqnarray}
Since, from $(1)$ we have $\stackrel{\ast}{\nabla}_Z X\in\Gamma(T_\lambda)~\mbox{for any}~ X,Z\in\Gamma(T_\lambda)$, then by (\ref{lm22}) we get $g(\nabla_Z X,Y)=0$. 
So, from relations (\ref{lm21}) and (\ref{lm22}), we obtain
\begin{eqnarray}
(\mu - \lambda)g(\nabla_X Y,Z)=-(\mu - \lambda)g(\nabla_Z X,Y)=0.
\end{eqnarray}
Hence, if $\lambda\neq\mu$, then $\nabla_X Y\bot T_\lambda$. \\
Similarly, we prove that $\nabla_Y X\bot T_\lambda~\mbox{if}~\lambda\neq\mu$. \\
\par Now we prove the following theorem which extends Cartan$^{\prime}$s fundamental formula on half-lightlike submanifold of Lorentzian manifolds with constant curvature.
\begin{thm}
\textsl{Let $(M,g,S(TM),S(TM^\bot))$ be an irrotational half-lightlike submanifold of an $(m+2)$-dimensional Lorenzian manifold $(\overline{M}(k),\overline{g})$ of constant curvature $k$. Assume that $E_0=\xi,E_1,...,E_{m-1}$ are eigenvectors of $\stackrel{\ast}{A}_\xi$ satisfying $\stackrel{\ast}{A}_\xi E_0=0~~\mbox{and} \stackrel{\ast}{A}_\xi E_i=\lambda_i E_i$ such that $\lambda_i$ is constant along $S (TM)$ for all $i$ and $\tau=0$ ($\{E_i\}_{i=1,...m-1}$ represents an orthonormal basis of $S(TM)$). Then, for every $i\in\{1,...,m-1\}$, we have:
\begin{eqnarray}
\sum_{\substack{j=1 \\ \lambda_j\neq \lambda_i}}^{m-1}\frac{k+\lambda_j g(A_N E_i,E_i)+\lambda_i g(A_N E_j,E_j)+g(A_L E_j,E_j)g(A_L E_i,E_i)-g(A_L E_i,E_j)^2}{\lambda_i - \lambda_j}=0
\end{eqnarray}
Moreover, if the screen is conformal with conformal factor $\varphi$, then for all $\in\{1,\cdots, m-1\}$,
\begin{eqnarray}
\sum_{\substack{j=1 \\ \lambda_j\neq \lambda_i}}^{m-1}\frac{k+2\varphi\lambda_i\lambda_j+g(A_L E_j,E_j)g(A_L E_i,E_i)-g(A_L E_i,E_j)^{2}}{\lambda_i - \lambda_j}=0
\end{eqnarray}
}
\end{thm}
\textbf{Proof}\\
From $(1)$ of proposition\ref{prop1}, and relations (\ref{shape01}) and (\ref{eqqconfcscren}), we have:
\begin{eqnarray}\label{theocartan1}
R(E_i,E_j)Ej = kE_i+\lambda_j A_N E_i +g(A_L E_j,E_j)A_L E_i-g(A_L E_i,E_j)A_L E_j
\end{eqnarray}
Also, we have:
\begin{eqnarray}
R(E_i,E_j)E_j = \nabla_{E_i}\nabla_{E_j} E_j-\nabla_{E_j}\nabla_{E_i} E_j-\nabla_{[E_i,E_j]}E_j
\end{eqnarray}
then, using (\ref{metric}) and lemma\ref{lm2} for $\lambda_i\neq\lambda_j$, we get:
\begin{eqnarray}
g(R(E_i,E_j)E_j,E_i)&=& g(\nabla_{E_i}\nabla_{E_j}E_j,E_i)-g(\nabla_{E_j}\nabla_{E_i} E_j,E_i)-g(\nabla_{[E_i,E_j]}E_j,E_i)
\cr&=& g(\nabla_{E_i}\nabla_{E_j} E_j,E_i)+(\nabla_{E_j}g)(\nabla_{E_i}E_j,E_i)-E_j.g(\nabla_{E_i}E_j,E_i)+g(\nabla_{E_i}E_j,\nabla_{E_j}E_i)
\cr& & -g(\nabla_{[E_i,E_j]}E_j,E_i)
\cr&=& g(\nabla_{E_i}\nabla_{E_j} E_j,E_i)+g(\nabla_{E_i}E_j,\nabla_{E_j}E_i)-g(\nabla_{[E_i,E_j]}E_j,E_i)
\cr&=& E_i.g(\nabla_{E_j} E_j,E_i)-g(\nabla_{E_j} E_j,\nabla_{E_i} E_i)-(\nabla_{E_i}g)(\nabla_{E_j} E_j,E_i)+g(\nabla_{E_i}E_j,\nabla_{E_j}E_i)
\cr& & -g(\nabla_{[E_i,E_j]}E_j,E_i)
\cr&=& -g(\nabla_{E_j} E_j,\nabla_{E_i} E_i)-B(E_i,E_i)g(\nabla_{E_j} E_j,N)+g(\nabla_{E_i} E_j,\nabla_{E_j} E_i)-g(\nabla_{[E_i,E_j]}E_j,E_i)
\cr&=& -g(\stackrel{\ast}{\nabla}_{E_j} E_j,\stackrel{\ast}{\nabla}_{E_i} E_i)-\lambda_i C(E_j,E_j)+g(\nabla_{E_i} E_j,\nabla_{E_j} E_i)-g(\nabla_{[E_i,E_j]}E_j,E_i)
\cr&=& -E_j.g(E_j,\stackrel{\ast}{\nabla}_{E_i} E_i)+g(E_j,\stackrel{\ast}{\nabla}_{E_j}\stackrel{\ast}{\nabla}_{E_i} E_i)-\lambda_i g(A_N E_j,E_j)+g(\nabla_{E_i} E_j,\nabla_{E_j} E_i)
\cr& & -g(\nabla_{[E_i,E_j]}E_j,E_i)
\cr&=& g(\nabla_{E_i} E_j,\nabla_{E_j} E_i)-g(\nabla_{[E_i,E_j]}E_j,E_i)-\lambda_i g(A_N E_j,E_j)
\label{theocartan2}
\end{eqnarray}
Otherwise, from (\ref{theocartan1}), we have:
\begin{eqnarray}\label{theocartan3}
g(R(E_i,E_i)E_j,E_i)&=& k+\lambda_j g(A_N E_i,E_i)+g(A_L E_j,E_j)g(A_L E_i,E_i)-g(A_L E_i,E_j)g(A_L E_j,E_i)
\end{eqnarray}
by (\ref{theocartan2}) and (\ref{theocartan3}), we have:
\begin{eqnarray}\label{theocartan4}
k+\lambda_j g(A_N E_i,E_i)+\lambda_i g(A_N E_j,E_j)+g(A_L E_j,E_j)g(A_L E_i,E_i)-g(A_L E_i,E_j)^{2}=
g(\nabla_{E_i} E_j,\nabla_{E_j} E_i)\nonumber \\-g(\nabla_{[E_i,E_j]}E_j,E_i)
\end{eqnarray}
Using $(v)$ of lemma\ref{lm1}, we get:
\begin{eqnarray}\label{theocartan5}
g((\nabla_{[E_i,E_j]}\stackrel{\ast}{A}_\xi)E_i,E_j)&=& (\lambda_i-\lambda_j)g(\nabla_{[E_i,E_j]}E_i,E_j)
\cr&=& (\lambda_j-\lambda_i)g(\nabla_{[E_i,E_j]}E_j,E_i)
\end{eqnarray}
then
\begin{eqnarray}\label{theocartan6}
g(\nabla_{[E_i,E_j]}E_j,E_i)=\frac{g((\nabla_{[E_i,E_j]}\stackrel{\ast}{A}_\xi)E_i,E_j)}{\lambda_j-\lambda_i}
\end{eqnarray}
Also, using $(i)$ and $(iv)$ of lemma\ref{lm1}, we have:
\begin{eqnarray}\label{theocartan7}
g((\nabla_{[E_i,E_j]}\stackrel{\ast}{A}_\xi)E_i,E_j)&=& g((\nabla_{E_i}\stackrel{\ast}{A}_\xi)E_j,[E_i,E_j])
\cr&=& g((\nabla_{E_i}\stackrel{\ast}{A}_\xi)E_j,\nabla_{E_i}E_j)-g((\nabla_{E_i}\stackrel{\ast}{A}_\xi)E_j,\nabla_{E_j}E_i)
\cr&=& g((\nabla_{E_j}\stackrel{\ast}{A}_\xi)E_i,\nabla_{E_i}E_j)-g((\nabla_{E_i}\stackrel{\ast}{A}_\xi)E_j,\nabla_{E_j}E_i)
\cr&=& g(\nabla_{E_j}(\stackrel{\ast}{A}_\xi E_i),\nabla_{E_i}E_j)-g(\stackrel{\ast}{A}_\xi\nabla_{E_j}E_i,\nabla_{E_i}E_j)-g(\nabla_{E_i}(\stackrel{\ast}{A}_\xi E_j),\nabla_{E_j}E_i)
\cr&  & +g(\stackrel{\ast}{A}_\xi\nabla_{E_i}E_j,\nabla_{E_j}E_i)
\cr&=& (\lambda_i-\lambda_j)g(\nabla_{E_i}E_j,\nabla_{E_j}E_i)
\end{eqnarray}
By (\ref{theocartan4}),(\ref{theocartan6}) and (\ref{theocartan7}), we get:
\begin{eqnarray}\label{theocartan8}
k+\lambda_j g(A_N E_i,E_i)+\lambda_i g(A_N E_j,E_j)+g(A_L E_j,E_j)g(A_L E_i,E_i)-g(A_L E_i,E_j)^{2}=2g(\nabla_{E_i}E_j,\nabla_{E_j}E_i)
\end{eqnarray}
Since $\nabla_{E_i}E_j=\sum_{s=1}^{m-1}g(\nabla_{E_i}E_j,E_s)E_s+\eta(\nabla_{E_i}E_j)\xi$, (\ref{theocartan8}) became
\begin{eqnarray}
k+\lambda_j g(A_N E_i,E_i)+\lambda_i g(A_N E_j,E_j)+g(A_L E_j,E_j)g(A_L E_i,E_i)-g(A_L E_i,E_j)^{2}=\nonumber\\ 2\sum_{s=1}^{m-1}g(\nabla_{E_i}E_j,E_s)g(\nabla_{E_j}E_i,E_s)
\end{eqnarray}
$(i)$ et $(v)$ of lemma\ref{lm1} give
\begin{eqnarray}
&&k+\lambda_j g(A_N E_i,E_i)+\lambda_i g(A_N E_j,E_j)+g(A_L E_j,E_j)g(A_L E_i,E_i)-g(A_L E_i,E_j)^{2}
\cr&=&2\sum_{\substack{s=1 \\ \lambda_s\neq \lambda_i,\lambda_j}}^{m-1}\frac{g((\nabla_{E_i}\stackrel{\ast}{A}_\xi)E_j,E_s)^2}{(\lambda_i-\lambda_s)(\lambda_j-\lambda_s)}
\end{eqnarray}
thus, we have
\begin{eqnarray}
\sum_{\substack{j=1 \\ \lambda_j\neq \lambda_i}}^{m-1}\frac{k+\lambda_j g(A_N E_i,E_i)+\lambda_i g(A_N E_j,E_j)+g(A_L E_j,E_j)g(A_L E_i,E_i)-g(A_L E_i,E_j)^{2}}{\lambda_i - \lambda_j}=\nonumber\\ \sum_{\substack{j=1 \\ \lambda_j\neq \lambda_i}}^{m-1}\frac{1}{\lambda_i -\lambda_j}2\sum_{\substack{s=1 \\ \lambda_s\neq \lambda_i,\lambda_j}}^{m-1}\frac{g((\nabla_{E_i}\stackrel{\ast}{A}_\xi)E_j,E_s)^2}{(\lambda_i-\lambda_s)(\lambda_j-\lambda_s)}\nonumber\\
= \sum_{\substack{s=1 \\ \lambda_s\neq \lambda_i}}^{m-1}-\frac{1}{\lambda_i -\lambda_s}2\sum_{\substack{j=1 \\ \lambda_j\neq \lambda_i,\lambda_s}}^{m-1}\frac{g((\nabla_{E_i}\stackrel{\ast}{A}_\xi)E_j,E_s)^2}{(\lambda_i-\lambda_j)(\lambda_s-\lambda_j)}\nonumber\\
= -\sum_{\substack{s=1 \\ \lambda_s\neq \lambda_i}}^{m-1}\frac{k+\lambda_j g(A_N E_i,E_i)+\lambda_i g(A_N E_j,E_j)+g(A_L E_j,E_j)g(A_L E_i,E_i)-g(A_L E_i,E_j)^{2}}{\lambda_i - \lambda_s}
\end{eqnarray}
Finally, we have
\begin{eqnarray}\label{theocartan9}
\sum_{\substack{j=1 \\ \lambda_j\neq \lambda_i}}^{m-1}\frac{k+\lambda_j g(A_N E_i,E_i)+\lambda_i g(A_N E_j,E_j)+g(A_L E_j,E_j)g(A_L E_i,E_i)-g(A_L E_i,E_j)^{2}}{\lambda_i-\lambda_j}=0
\end{eqnarray}
\section{Application}
A vector field $X$ on a semi-Riemannian manifold $(\overline{M},\overline{g})$ is said to be conformal vector field if there exist a smooth function $\sigma$ on $\overline{M}$ called potential function such that $\mathcal{L}_{X}\overline{g}=2\sigma\overline{g}$, , where the symbole $\mathcal{L}_{X}$ denote the Lie derivative with respect to the vector field $X$, that is,
\begin{equation}\label{killeq1}
(\mathcal{L}_{X}\overline{g})(Y,Z)= \overline{g}(\overline{\nabla}_{Y}X,Z)+\overline{g}(Z,\overline{\nabla}_{Z}X),
\end{equation}
for any $X,Y,Z\in\Gamma(T\overline{M})$. In particular, when a potential function $\sigma$ associated to a conformal vector field $X$ is identically null, we said that $X$ is a Killing vector field. A distribution $\mathcal{D}$ on $\overline{M}$ is said to be a conformal distribution if each vector field belonging to $\mathcal{D}$ is a conformal vector field.

 If the co-screen distribution $S(TM^{\bot})$ is a conformal distribution, using (\ref{GW3}) and (\ref{shape02}) we have 
\begin{eqnarray}
\overline{g}(\overline{\nabla}_{X}L,Y)=g(A_{L}X,Y)+\phi(X)\eta(Y)=-D(X,Y).
\end{eqnarray}
Therefore, we obtain
\begin{eqnarray}
(\mathcal{L}_{L}\overline{g})(X,Y)=-2D(X,Y).
\end{eqnarray}
A vector field $X$ on a semi-Riemannian manifold $(\overline{M},\overline{g})$ is said to be a Killing vector field if $\mathcal{L}_{X}\overline{g}=0$.

\begin{dfn}
A distribution $D$ on a semi-riemannian manifold $(\overline{M},\overline{g})$ is said to be a Killing distribution if each vector field belonging to $D$ is a killing vector field.
\end{dfn}
For a half-lightlike submanifold $(M,S(TM),S(TM^{\bot}))$ of a semi-Riemannian manifold $(\overline{M},\overline{g})$, it is easy to see that the co-screen distribution $S(TM^{\bot})$ is a Killing distribution if and only if $D(X,Y)=0$ for any $X,Y\in \Gamma(TM)$.

The if $S(TM^{\bot})$ is a Killing distribution we have:
\begin{eqnarray}\label{killeq5}
D(X,Y)=0,~~\phi(X)=0,~~A_{L}X=\rho(X)\xi,
\end{eqnarray}
for all $X,Y\in\Gamma(TM)$.
\begin{prop}\label{propconfcosreen}
Let $(M,g,S(TM),S(TM^\bot))$ be a half-lightlike submanifold of a semi-Riemannian manifold $(\overline{M},\overline{g})$. Then, the co-screen distribution $S(TM^{\bot})$ is a conformal distribution if and only if $D(X,Y)=-\sigma g(X,Y)$ for any $X,Y\in\Gamma(TM)$. 
\end{prop}
 Let $(M,g,S(TM),S(TM^\bot))$ be a half-lightlike submanifold of a semi-Riemannian manifold $(\overline{M}(k),\overline{g})$ of constant curvature $k$, with a conformal co-screen distribution. Then, by proposition\ref{propconfcosreen} and relations (\ref{shape02})and (\ref{secondformfondprop}), we have:
\begin{eqnarray}\label{eqqconfcscren}
D(X,Y)=-\sigma g(X,Y),~ \varphi(X)=0,~ A_{L}X=-\sigma PX+\rho(X)\xi,~ \mbox{for all}~ X,Y\in\Gamma(TM).
\end{eqnarray}
\par A semi-Riemannian manifold $(\overline{M},\overline{g})$ is said to be equipped with a non-trivial closed conformal vector field $X$, if for a certain smooth function $\sigma\in\mathcal{C}^{\infty}(\overline{M})$, we have 
\begin{eqnarray}
\overline{\nabla}_{Y}X=\sigma Y~\mbox{for any}~~Y\in\Gamma(T\overline{M}).
\end{eqnarray}
Where $\overline{\nabla}$ is the Levi-Civita connection associated to the metric $\overline{g}$ on $\overline{M}$.
\par It is obvious that a closed conformal vector field is a conformal vector field, but the opposite way is wrong. From (\ref{GW2}) and (\ref{eqqconfcscren}), we have the following
\begin{prop}
Let $(M,g,S(TM),S(TM^\bot))$ be a half-lightlike submanifold of a semi-Riemannian manifold $(\overline{M},\overline{g})$, with a conformal co-screen distribution $S(TM^{\bot})=\{L\}$. Then, the vector field $L$ is a closed conformal vector field if and only if $\rho(X)=-\sigma\eta(X)~\mbox{for any}~X\in\Gamma(TM)$.  
\end{prop} 

\begin{thm}\label{thm1conf}
Let $(M,g,S(TM),S(TM^\bot))$ be a screen homothetic half-lightlike submanifold of a semi-Riemannian manifold $(\overline{M}(k),\overline{g})$ of constant curvature $k$, with a conformal co-screen distribution. Then, we obtain
\begin{eqnarray}
2\varphi\tau(\xi)B=\{-k-\sigma^{2}\}g.
\end{eqnarray}
\end{thm}
\textbf{Proof}
Let $X,Y,Z\in\Gamma(TM)$. As $\overline{M}$ has a constant curvature $k$, we have $\overline{g}(\overline{R}(X,Y)Z,\xi)=0$. Then,  considering the relation (\ref{GausCoda1}), we get:
\begin{eqnarray}\label{proofprop1confeq1}
(\nabla_{X}B)(Y,Z)-(\nabla_{Y}B)(X,Z)=B(X,Z)\tau(Y)-B(Y,Z)\tau(X).
\end{eqnarray}
Also, using relations (\ref{GausCoda1}) and (\ref{GausCoda4}), we have:
\begin{eqnarray}
\overline{g}(R(X,Y)Z,N)=k\{g(Y,Z)\eta(X)-g(X,Z)\eta(Y)\}+\sigma \rho(Y)g(X,Z)-\sigma\rho(X)g(Y,Z),\label{proofprop1confeq2}
\end{eqnarray}
and
\begin{eqnarray}
\overline{g}(R(X,Y)PZ,N)&=&(\nabla_{X}C)(Y,PZ)-(\nabla_{Y}C)(X,PZ)+C(X,PZ)\tau(Y)-C(Y,PZ)\tau(X)
\cr&=& \varphi(\nabla_{X}B)(Y,PZ)-\varphi(\nabla_{Y}B)(X,PZ)+\varphi B(X,Z)\tau(Y)-\varphi B(Y,Z)\tau(X)
\cr&\overset{\ref{proofprop1confeq1}}{=}& 2\varphi\{B(X,PZ)\tau(Y)-B(Y,PZ)\tau(X)\}
\cr&=&2\varphi\{B(X,Z)\tau(Y)-B(Y,Z)\tau(X)\}\label{proofprop1confeq3}
\end{eqnarray}
Since $\overline{g}(R(X,Y)\xi,N)=0$, then $\overline{g}(R(X,Y)Z,N)=\overline{g}(R(X,Y)PZ,N)$. 
Therefore, from relations (\ref{proofprop1confeq2}) and (\ref{proofprop1confeq3}), we have:
\begin{eqnarray}
2\varphi\{B(X,Z)\tau(Y)-B(Y,Z)\tau(X)\}=k\{g(Y,Z)\eta(X)-g(X,Z)\eta(Y)\}+\sigma \rho(Y)g(X,Z)\nonumber\\
-\sigma\rho(X)g(Y,Z).
\end{eqnarray}
Replacing $Y$ by $\xi$ in the last relation, we obtain this theorem.
\begin{rque}\label{rque2}
\textsl{For a screen homothetic half-lightlike submanifold $(M,g,S(TM),S(TM^\bot))$ of a semi-Riemannian manifold $(\overline{M}(k),\overline{g})$, with a conformal co-screen distribution, the theorem\ref{thm1conf} allow us to show that, if $\tau=0$ then $k=-\sigma^{2}$, and if $\tau(\xi)\neq 0$ then $M$ is a totally umbilical submanifold of $\overline{M}(k)$.}
\end{rque}

If the ambient manifold $\overline{M}$ is a semi-Riemannian manifold $(\overline{M}(k),\overline{g})$ of constante curvature 
$k$, then $\overline{R}(\xi,Y)X=k\overline{g}(X,Y)\xi$; $\overline{R}(L,X)Y=k\overline{g}(X,Y)L$ and 
$\overline{R}ic(X,Y)=(m+1)k\overline{g}(X,Y)$. Thus, we obtain
\begin{eqnarray}\label{ricdengeneq9}
Ric(X,Y)=(m-1)kg(X,Y)+B(X,Y)tr(A_N)+D(X,Y)tr(A_L)-g(A_{N}X,\stackrel{\ast}{A}_{\xi}Y)\nonumber\\-g(A_{L}X,A_{L}Y)
+\rho(X)\varphi(Y)
\end{eqnarray}
A screen conformal half-lightlike submanifold of a semi-Riemannian manifold $(\overline{M}(k),\overline{g})$ of constante curvature $k$, with a conformal co-screen distribution, have a symmetric Ricci tensor $Ric$, by (\ref{ricdengeneq9}).
\begin{thm}\label{thmRicconfccsd}
\textsl{Let $(M,g,S(TM),S(TM^\bot))$ be a screen conformal half-lightlike submanifold of a semi-Riemannian manifold 
$(\overline{M}(k),\overline{g})$  of constant curvature $k$, with a conformal co-screen distribution. Then, $M$ admits an induced symmetric tensor.}
\end{thm}

Now, we prove the following theorem,

\begin{thm}\label{thmappli1}
\textsl{Let $(M,g,S(TM)),S(TM^\bot))$ be a screen homothetic half-lightlike submanifold of an $(m+2)$-dimensional Lorenzian manifold $(\overline{M}(k),\overline{g})$ of non-positive constant curvature $k$, 
 with a conformal co-screen distribution $S(TM^{\bot})$, whose screen principal curvatures along the screen distribution $S(TM)$ are constant and at most two are distinct. If $M$ has two distinct screen  principal curvatures, then one of then must be zero.}
\end{thm}
\textbf{Proof}\\
If $M$ has two distinct screen  principal curvatures $\alpha~\mbox{and}~\beta$, then it exists $p\in\{1,\cdots,m-1\}$ so that 
\begin{eqnarray}
\lambda_1=\lambda_2=\cdots=\lambda_p=\alpha~\mbox{and}~\lambda_{p+1}=\lambda_{p+2}=\cdots=\lambda_{m-1}=\beta
\end{eqnarray} 
Since $M$ is screen homothetic with a conformal co-screen distribution, then using relations (\ref{conf1}) and (\ref{eqqconfcscren}), the relation (\ref{theocartan8}) became 
\begin{eqnarray}
k+\sigma^{2}+2\varphi\alpha\beta=2g(\nabla_{E_\alpha}E_\beta,\nabla_{E_\beta}E_\alpha).
\end{eqnarray}
Then, according to the remark\ref{rque2}, we get:
\begin{eqnarray}
\varphi\alpha\beta=g(\nabla_{E_\alpha}E_\beta,\nabla_{E_\beta}E_\alpha)
\end{eqnarray}
From (2) of lemme\ref{lm2}, $g(\nabla_{E_\alpha}E_\beta,\nabla_{E_\beta}E_\alpha)=0$. therefore, previous relation give us $\varphi\alpha\beta=0$. Since $\varphi\neq 0~\mbox{and}~\alpha\neq\beta$, then  $\alpha=0~\mbox{and}~\beta\neq 0$ or $\alpha\neq 0~\mbox{and}~\beta= 0$.
\begin{exple}
Let $M$ be an immersed submanifold in $\mathbb{R}_{1}^{5}$ defined by:
\begin{eqnarray}
M &\rightarrow &\mathbb{R}_{1}^{5}\nonumber\\
(v_1,v_2,v_3)&\mapsto & (x_1,x_2,x_3,x_4,x_5)
\end{eqnarray}
where~~
$\begin{cases}
x_1=v_1\\
x_2=v_2\\
x_3=v_3\\
x_4=\sqrt{v_{1}^{2}-v_{2}^{2}}\\
x_5=\sqrt{1+v_{3}^{2}}
\end{cases}$
with $v_1>v_2>0$
\vspace{0.5cm}
\par the tangent bundle $TM$ is  $TM=span\{V_1,V_2,V_3\}$, where
\begin{eqnarray}
V_1=\frac{1}{x_2}\left(x_4\frac{\partial}{\partial x_1}+x_1\frac{\partial}{\partial x_4}\right)\\
V_2=\frac{1}{x_1}\left(x_4\frac{\partial}{\partial x_2}-x_2\frac{\partial}{\partial x_4}\right)\\
V_3=\frac{1}{\sqrt{1+2x_{3}^{2}}}\left(x_5\frac{\partial}{\partial x_3}+x_3\frac{\partial}{\partial x_5}\right)
\end{eqnarray}
and the  normal bundle is $TM^{\bot}=span\{H_{1},H_{2}\}$, where
\begin{eqnarray}
H_{1}=x_{1}\frac{\partial}{\partial x_1}+x_{2}\frac{\partial}{\partial x_2}+x_{4}\frac{\partial}{\partial x_4}\\
H_{2}=\frac{1}{\sqrt{1+2x_{3}^{2}}}\left(x_{3}\frac{\partial}{\partial x_3}+x_{5}\frac{\partial}{\partial x_5}\right)
\end{eqnarray}
\begin{eqnarray}
\overline{g}(H_{1},H_{1})=0~~\mbox{and}~~\overline{g}(H_{2},H_{2})=1,
\end{eqnarray}
then $H_{1}\in Rad (TM)$and $H_2\in  S(TM^{\bot})$. Thus $ Rad( TM)$ and $S(TM^{\bot})$ are of rank  1. 
Remark that $H_1=\frac{x_{1}x_{2}}{x_4}(V_{1}+V_{2})$. 
\par Take $ Rad (TM)=span\{\xi\}$ where $\xi=\frac{1}{x_{1}}H_{1}=\frac{1}{x_{1}}\{x_1\frac{\partial}{\partial x_1}+x_{2}\frac{\partial}{\partial x_2}+x_{4}\frac{\partial}{\partial x_4}\}$, and $S(TM^{\bot})=span\{L=H_{2}\}$ 
 
the null transversal vector field  $N$ is:
\begin{eqnarray}
N = -\frac{1}{ 2x_1}\left\{x_1\frac{\partial}{\partial x_1}-x_{2}\frac{\partial}{\partial x_2}-x_{4}\frac{\partial}{\partial x_4}\right\},
\end{eqnarray}
The null transversal bundle $ltr(TM)$ and the  screen distribution $S(TM)$ are:
\begin{eqnarray}
ltr(TM)=span\{N\}~\mbox{and}~S(TM)=span\{V_2,V_3\}.
\end{eqnarray}
By direct computation, we get
\begin{eqnarray}
\overline{\nabla}_{V_2}\xi=\frac{1}{ x_{1}}V_2,~~\overline{\nabla}_{V_3}\xi=0,~~\overline{\nabla}_{\xi}\xi=0
\end{eqnarray}
Thus, 
\begin{eqnarray}
\stackrel{\ast}{A}_{\xi}V_2=-\frac{1}{ x_{1}}V_2,~~\stackrel{\ast}{A}_{\xi}V_3=0,~~\stackrel{\ast}{A}_{\xi}\xi=0,~~\tau=0.
\end{eqnarray}
Then, $M$ is \textbf{irrotational} and has \textbf{two distinct screen principal curvatures} $\lambda=-\frac{1}{x_1}~\mbox{and}~\mu=0$. On the other hand, it is easy to see that:
\begin{eqnarray}
\overline{\nabla}_{V_2}N= \frac{1}{2x_{1}}V_2,~~\overline{\nabla}_{V_3}N=0,~~\overline{\nabla}_{\xi}N=0
\end{eqnarray}
then,
\begin{eqnarray}
A_{N}V_2=-\frac{1}{2x_{1}}V_2,~~A_{N}V_3=0,~~A_{N}\xi=0,~~\rho=0.
\end{eqnarray}
it follow
\begin{eqnarray}
A_{N}V_2=\frac{1}{2}\stackrel{\ast}{A}_{\xi}V_2,~~A_{N}V_3=\stackrel{\ast}{A}_{\xi}V_3=0,~~A_{N}\xi=\stackrel{\ast}{A}_{\xi}\xi=0,
\end{eqnarray}
which proves that $M$ is irrotational screen homothetic half-lightlike submanifold, with  conformal factor $\frac{1}{2}$ and the 1-forme $\tau=0$.
\end{exple}
\begin{exple}
Let $M$ be a submanifold in $\mathbb{R}_{1}^{6}$ given by:

\begin{equation}
x_3=x_5,\;\;\; x_6=\sqrt{x_{1}^{2}-x_{2}^{2}}\;\;\;\mbox{with}~~x_1>x_2>0.
\end{equation}

 Then, we have:
\begin{eqnarray}
Rad (TM)&=&span\left\{\xi=\frac{1}{x_{1}}\left(x_1\frac{\partial}{\partial x_1}+x_{2}\frac{\partial}{\partial x_2}+x_{6}\frac{\partial}{\partial x_6}\right)\right\}
\cr S(TM^\bot)&=&span\left\{L=-\frac{1}{\sqrt{2}}\left(-\frac{\partial}{\partial x_3}+\frac{\partial}{\partial x_5}\right)\right\}
\cr ltr(TM)&=&span\left\{-\frac{1}{ 2x_1}\left(x_1\frac{\partial}{\partial x_1}-x_{2}\frac{\partial}{\partial x_2}-x_{6}\frac{\partial}{\partial x_6}\right)\right\}
\cr S(TM)&=&span\left\{V_1, V_2, V_3\right\}
\end{eqnarray}
 where
\begin{eqnarray}
V_1&=&\frac{1}{x_1}\left(x_6\frac{\partial}{\partial x_2}-x_2\frac{\partial}{\partial x_6}\right)
\cr V_2&=&\frac{1}{\sqrt{2}}\left(\frac{\partial}{\partial x_3}+\frac{\partial}{\partial x_5}\right)
\cr V_3&=&\frac{\partial}{\partial x_4}
\end{eqnarray}
By straightforward calculation, we get
\begin{eqnarray}
\overline{\nabla}_{V_1}\xi= \frac{1}{ x_{1}}V_1,~~\overline{\nabla}_{V_2}\xi=0;~~\overline{\nabla}_{V_3}\xi=0;~~\overline{\nabla}_{\xi}\xi=0
\end{eqnarray}
Thus, $M$ is \textbf{irrotational} and 
\begin{eqnarray}\label{eq1exemple2}
\stackrel{\ast}{A}_{\xi}V_1=-\frac{1}{ x_{1}}V_1,~~\stackrel{\ast}{A}_{\xi}V_2=\stackrel{\ast}{A}_{\xi}V_3=\stackrel{\ast}{A}_{\xi}\xi=0,~~\tau=0.
\end{eqnarray}
Then, $M$ has \textbf{two distinct screen principal curvatures} $\lambda=-\frac{1}{x_1}~\mbox{and}~\mu=0$. On the other hand, we have:
\begin{eqnarray}
\overline{\nabla}_{V_1}N= \frac{1}{2x_{1}}V_1,~~\overline{\nabla}_{V_2}N=0,~~\overline{\nabla}_{V_3}N=0,~~\overline{\nabla}_{\xi}N =0
\end{eqnarray}
it follow that
\begin{eqnarray}\label{eq2exemple2}
A_{N}V_1=-\frac{1}{2x_{1}}V_1,~~A_{N}V_2=A_{N}V_3=A_{N}\xi=0,~\mbox{and}~\rho=0.
\end{eqnarray}
By (\ref{eq1exemple2}) and (\ref{eq2exemple2}) $M$ is a screen homothetic half-lightlike submanifold, with  conformal factor $\frac{1}{2}$.
\par Otherwise, we have:
\begin{eqnarray}
\overline{\nabla}_{\xi}L=\overline{\nabla}_{V_2}L=\overline{\nabla}_{V_3}L=\overline{\nabla}_{V_4}L=0,
\end{eqnarray}
and then,
\begin{eqnarray}
A_{L}\equiv 0~\mbox{on}~TM.
\end{eqnarray}
Consequently,
\begin{eqnarray}
D(X,Y)=g(A_{L}X,Y)=0~\mbox{for any}~X,Y\in\Gamma(TM),
\end{eqnarray}
which show that the co-screen distribution $S(TM)=span\{L\}$ is a Killing distribution. 

\end{exple}
\begin{exple}
Let $M$ be an immersed submanifold in $\mathbb{R}_{1}^{6}$ defined by:
\begin{eqnarray}
M &\rightarrow &\mathbb{R}_{1}^{6}\nonumber\\
(v_1,v_2,v_3)&\mapsto & (x_1,x_2,x_3,x_4,x_5,x_6)
\end{eqnarray}
where~~
$\begin{cases}
x_1=v_1\\
x_2=v_2\\
x_3=\frac{v_1}{\sqrt{2}}\sin{v_3}\\
x_4=\frac{v_1}{\sqrt{2}}\cos{v_3}\\
x_5=\frac{v_1}{\sqrt{2}}\sin{v_4}\\
x_6=\frac{v_1}{\sqrt{2}}\cos{v_4}\\
\end{cases}$
with $v_1\neq 0~~\mbox{and}~~v_3,v_4\in\mathbb{R}\backslash\{\frac{\pi}{2}+k\pi\}$
\vspace{0.5cm}
\par the tangent bundle $TM$ is  $TM=span\{V_1,V_2,V_3,V_4\}$, where
\begin{eqnarray}
V_1&=&\frac{1}{x_1}\left\{x_1\frac{\partial}{\partial x_1}+x_3\frac{\partial}{\partial x_3}+x_4\frac{\partial}{\partial x_4}+x_5\frac{\partial}{\partial x_5}+x_6\frac{\partial}{\partial x_6}\right\}
\cr V_2&=&\frac{\partial}{\partial x_2}\\
\cr V_3&=&\frac{\sqrt{2}}{x_1}\left\{x_4\frac{\partial}{\partial x_3}-x_3\frac{\partial}{\partial x_4}\right\}\\
\cr V_4&=&\frac{\sqrt{2}}{x_1}\left\{ x_6\frac{\partial}{\partial x_5}-x_5\frac{\partial}{\partial x_6}\right\}
\end{eqnarray}
and the  normal bundle is $TM^{\bot}=span\{H_{1},H_{2}\}$, where
\begin{eqnarray}
H_{1}=\frac{1}{2x_4}\left\{x_{1}\frac{\partial}{\partial x_1}+2x_{3}\frac{\partial}{\partial x_3}+2x_{4}\frac{\partial}{\partial x_4}\right\}\\
H_{2}=\frac{1}{2x_6}\left\{x_{1}\frac{\partial}{\partial x_1}+2x_{5}\frac{\partial}{\partial x_5}+2x_{6}\frac{\partial}{\partial x_6}\right\}
\end{eqnarray}
Remark that $\frac{x_4}{x_1}H_1+\frac{x_6}{x_1}H_2=V_1$. 
\par Take $ Rad (TM)=span\{\xi=V_1\}$, $S(TM)=span\{V_2,V_3,V_4\}$ and $S(TM^{\bot})=span\{L=H_1\}$.
the null transversal vector field N is:
\begin{eqnarray}
N= -\frac{1}{ x_1}\left\{x_1\frac{\partial}{\partial x_1}+x_3\frac{\partial}{\partial x_3}+x_4\frac{\partial}{\partial x_4}-x_5\frac{\partial}{\partial x_5}-x_6\frac{\partial}{\partial x_6}\right\}
\end{eqnarray}
The null transversal bundle $ltr(TM)$ is ltr(TM)=span\{N\}.

By direct computation we get:
\begin{equation}
\overline{\nabla}_{\xi}\xi=0,\;\;\;\overline{\nabla}_{V_2}\xi=0,\;\;\; \overline{\nabla}_{V_3}\xi=\frac{1}{x_1}V_3,\;\;\; \overline{\nabla}_{V_4}\xi=\frac{1}{x_1}V_4
\end{equation}

Thus, we deduce that $M$ is \textbf{irrotational} and:
\begin{eqnarray}
\stackrel{\ast}{A}_{\xi}\xi=0,~~\stackrel{\ast}{A}_{\xi}V_2=0,~~\stackrel{\ast}{A}_{\xi}V_3=-\frac{1}{ x_1}V_3,~~\stackrel{\ast}{A}_{\xi}V_4=-\frac{1}{ x_1}V_4~~\mbox{and}~~\tau=0.
\end{eqnarray}
Then, $M$ has \textbf{two distinct screen principal curvatures} $\lambda=-\frac{1}{x_1}~\mbox{and}~\mu=0$. On the other hand, we have:

\begin{equation}
\overline{\nabla}_{\xi}N=0,\;\;\;\overline{\nabla}_{V_2}N=0,\;\;\; \overline{\nabla}_{V_3}N=-\frac{1}{x_1}V_3,\;\;\; \overline{\nabla}_{V_4}N=\frac{1}{x_1}V_4
\end{equation}

Then,
\begin{eqnarray}
A_{N}\xi=0,~~ A_{N}V_2=0,~~A_{N}V_3=\frac{1}{ x_{1}}V_3,~~ A_{N}V_4=-\frac{1}{ x_{1}}V_4~~\mbox{and}~~\rho=0.
\end{eqnarray}
It follow that,
\begin{eqnarray}
A_{N}\xi=\stackrel{\ast}{A}_{\xi}\xi=0,~~ A_{N}V_2=\stackrel{\ast}{A}_{\xi}V_2=0,~~A_{N}V_3=-\stackrel{\ast}{A}_{\xi}V_3,~~A_{N}V_4=\stackrel{\ast}{A}_{\xi}V_4.
\end{eqnarray}
So, $M$ is not a screen homothetic half-lightlike submanifold.
\par Otherwise, we have:
\begin{eqnarray}
D(V_2,V_2)=\overline{g}(\overline{\nabla}_{V_2}V_2,L)=0,
\end{eqnarray}
\begin{eqnarray}
D(V_2,V_3)=\overline{g}(\overline{\nabla}_{V_2}V_3,L)=0,
\end{eqnarray}
\begin{eqnarray}
D(V_2,V_4)=\overline{g}(\overline{\nabla}_{V_2}V_4,L)=0,
\end{eqnarray}
\begin{eqnarray}
D(V_3,V_3)=\overline{g}(\overline{\nabla}_{V_3}V_3,L)=-\frac{\sqrt{2}}{x_{1}x_4}g(V_3,V_3),
\end{eqnarray}
\begin{eqnarray}
D(V_3,V_4)=\overline{g}(\overline{\nabla}_{V_3}V_4,L)=0,
\end{eqnarray}
\begin{eqnarray}
D(V_4,V_4)=\overline{g}(\overline{\nabla}_{V_4}V_4,L)=-\frac{\sqrt{2}}{x_{1}}g(V,L)=0.
\end{eqnarray}
Which  show that the co-screen $S(TM^\bot)=span\{L=H_1\}$ is not conformal.
\par This shows that a screen conformal and co-screen conformal distribution are just a sufficient condition for a half-lightlike submanifold with two distinct screen principal curvatures then one must be zero.
\end{exple}
 In the sequel, we consider a screen homothetic half-lightlike submanifold $M$ of an $(m+2)$-dimensional Lorenzian manifold $(\overline{M}(k),\overline{g})$ of non-positive constant curvature $k$, 
with a conformal co-screen distribution $S(TM^{\bot})$, whose screen principal curvatures are constant along $S(TM)$. We assume that $M$ has exactly two distinct screen principal curvatures. Then, by theorem\ref{thmappli1}, one of them must be zero.
We denote by $\lambda$ the non-zero screen principal curvature and $r$ the multiplicity of $\lambda$. The sets, 
\begin{eqnarray}
T_\lambda = \{X\in\Gamma(S(TM)):~~\stackrel{\ast}{A}_\xi X=\lambda X \}
\end{eqnarray}
\begin{eqnarray}
T_0 = \{X\in\Gamma(S(TM)):~~\stackrel{\ast}{A}_\xi X=0 \}
\end{eqnarray}
define the distributions of dimensions $r$ and dimension $m-r$, respectively.
\par From (\cite{K.S}) and the remark\ref{rque1}, it is obvious that if $M$ is a screen conformal half-lightlike submanifold of a Lorentzian manifold $(\overline{M}(k),\overline{g})$, with a conformal co-screen distribution, then the screen distribution $S(TM)$ is Riemannian and integrable. By the theorem\ref{thmRicconfccsd}, we have that the induced Ricci tensor on $M$ is symmetric. Moreover, a screen conformal half-lightlike submanifold is locally a product $C\times M'$ where $C$ is a null curve, $M'$ is an integral manifold of $S(TM)$ (\cite{D.H.J1}). We have the following local decomposition.
\begin{thm}
\textsl{Let $(M,g,S(TM),S(TM^\bot))$ be a screen homothetic half-lightlike  submanifold of an $(m+2)$-dimensional Lorenzian manifold $(\overline{M}(k),\overline{g})$ of non-positive constant curvature $k$,
with a conformal co-screen distribution $S(TM^{\bot})$. If the principal curvatures of $M$ are constant along the screen distribution $S(TM)$ and exactly two of them are  distinct, then $M$ is locally a lightlike triple product manifold $M=C\times (M'=M_\lambda\times M_0)$, where C is a null curve, $M'$ is an integral manifold of $S(TM)$, $M_\lambda~\mbox{and}~M_0$ are leaves of some distributions of $M$ such that $M_\lambda$ is an $r$-dimensional totally geodesic Riemannian manifold of curvature $2\varphi\lambda^{2}$ and $M_0$ is an $(m-r)$-dimensional totally geodesic Euclidien espace.
}
\end{thm}
\textbf{Proof}\\
$S(TM)$ being Riemannian and integrable, it is know that its leaf $M'$ is Riemannian. Since $M$ has exactly two distinct screen principal curvatures, then one of them must be zero (theorem\ref{thmappli1}). 
 Then, with previous notations, we have 
 \begin{eqnarray}
S(TM)=T_\lambda\underset{\bot}{\oplus}T_0.
\end{eqnarray}
Now, we prove that $T_\lambda~\mbox{et}~T_0$ are integrables.
Let $X,Y\in\Gamma(T_\lambda)$, we have:
\begin{eqnarray}
\stackrel{\ast}{A}_\xi [X,Y]&=& \stackrel{\ast}{A}_\xi\nabla_X Y -\stackrel{\ast}{A}_\xi \nabla_Y X
\cr&=& \nabla_X(\stackrel{\ast}{A}_\xi Y)-(\nabla_X\stackrel{\ast}{A}_\xi)Y -\nabla_Y(\stackrel{\ast}{A}_\xi X)+ (\nabla_Y\stackrel{\ast}{A}_\xi)X
\cr&=& \nabla_X(\stackrel{\ast}{A}_\xi Y)-\nabla_Y(\stackrel{\ast}{A}_\xi X)
\cr&=& \lambda\nabla_X Y -\lambda \nabla_Y X
\cr&=& \lambda [X,Y],
\end{eqnarray}
that is $T_\lambda$ is involutive, then integrable.
 Similarly $T_0$ is integrable.\\
Using the point (2) of the lemma\ref{lm2}, we get that $T_\lambda~\mbox{and}~T_0$ are parallel along their normals in $S (TM)$. \\
By the decomposition theorem of de Rham(\cite{G.R}), we have $M'=M_\lambda\times M_0$; where $M_\lambda~\mbox{and}~M_0$ are some leaves of $T_\lambda~~\mbox{and}~~T_0$ respectively. Thus $M$ is locally a
product $C\times M'=C\times M_\lambda\times M_0$.\\
Remark that $M_\lambda$ is totally geodesique $\Leftrightarrow~ \stackrel{\ast}{\nabla}_X Y\in\Gamma(T_\lambda)~\mbox{for all}~ X,Y\in\Gamma(T_\lambda)$. Then, by (1) of the  lemma\ref{lm2},  $M_\lambda$ is totally geodesique. Similarly, we prove that $M_0$ is totally geodesique in $S(TM)$. \\

Consider the frame field of eigenvectors $\{E_1, E_2,\cdots, E_r\}$ of  $\stackrel{\ast}{A}_\xi$ such that $\{E_i\}_{i=1,\cdots,r}$ is an orthonormal frame field of $T_\lambda$, then using (\ref{GausCoda4}) et (\ref{theocartan1}) we have
\begin{eqnarray}
g(R(E_i,E_j)E_j,E_i)=\varphi\lambda^2=g(\stackrel{\ast}{R}(E_i,E_j)E_j,E_i)-\varphi\lambda^2
\end{eqnarray}
then
\begin{eqnarray}
g(\stackrel{\ast}{R}(E_i,E_j)E_j,E_i)=2\varphi\lambda^2
\end{eqnarray}
Thus the sectional curvature $K^\lambda$ of the leaf $M_\lambda$ of $T_\lambda$ is given by
\begin{eqnarray}
K^\lambda (E_i,E_j)=\frac{g(\stackrel{\ast}{R}(E_i,E_j)E_j,E_i)}{g(E_i,E_i)g(E_j,E_j)-g(E_i,E_j)^{2}}=2\varphi\lambda^2.
\end{eqnarray}
By the same way, we can see that the sectional curvature of the leaf $M_0$ of $T_0$ is 
\begin{eqnarray}
K^{0}=0. 
\end{eqnarray}
This completes the proof.

\par \textsl{Next, we say that a half-lightlike submanifold $M$ is $B$ (resp. $D$)-umbilical if on each coordinate neighborhood $\mathcal{U}\subset M$, there exists smooth functions $H_1$ (resp. $H_2$) such that
$B(X,Y) = H_1 g(X,Y)$~~(resp. $D(X,Y)=H_2 g(X,Y)$) for any $X,Y\in\Gamma(TM\vert_\mathcal{U})$. In particular, if $H_1=0$ (resp. $H_2=0$), then $M$ is said to be totally $B$ (resp. $D$)-geodesic.}
\par \textsl{Note that, $M$ is totally umbilical ( resp. totally geodesic ) if and only if $M$ are both $B$ and $D$-totally umbilical (resp. totally geodesic ).}
 \par It is easy to see that, for a half-lightlike submanifold  $M$ with a conformal co-screen distribution :
\begin{enumerate}
\item If the screen principal curvatures are all identical and non-zero then $M$ is totally umbilical,
\item If the screen principal curvatures are all identical all zero then $M$ is totally geodesic.
\end{enumerate}   . 
\begin{rque}
\textsl{ We assume that $M$ has at most two distinct screen principal curvatures, thus, if
the screen principal curvatures are all identical, $M$ is either totally geodesic or totally umbilical and
if the two screen principal curvatures are distinct, then $M=C\times M_\lambda\times M_\mu$}
\end{rque}
Thus we have the following classification theorem.
\begin{thm}
\textsl{Let $(M,g,S(TM),S(TM^\bot))$ be a screen homothetic half-lightlike  submanifold of an $(m+2)$-dimensional Lorenzian manifold $(\overline{M}(k),\overline{g})$ of non-positive constant curvature $k$,
with a conformal co-screen distribution $S(TM^{\bot})$. If the screen principal curvatures of $M$ are all constant along the screen distribution
$S(TM)$ such that at most two of them are distinct. Then we have one of the following:
\begin{enumerate}
\item $M$ is either totally geodesic or totally umbilical;
\item $M$ is locally a lightlike triple product manifold  $C\times M_\lambda\times M_0$, where $C$ is a null curve, $M_\lambda~\mbox{and}~M_0$ are leaves of some distributions of $M$  such that $M_\lambda$ is a totally geodesic Riemannian manifold of curvature $2\varphi\lambda^2$ and $M_0$ is an $(m-r)$-dimensional totally geodesic Euclidien espace. 
\end{enumerate}
}
\end{thm}

\end{document}